\documentclass[10pt]{article}
\usepackage{cite}
\usepackage{url}
\usepackage{ifthen}
\usepackage{multicol}
\usepackage{amsmath}
\usepackage{latexsym,amssymb,amsmath}
\usepackage[utf8]{inputenc}
\usepackage{graphicx}
\usepackage{color}

\newtheorem{theorem}{Theorem}
\newtheorem{corollary}{Corollary}
\newtheorem{lemma}{Lemma}
\newtheorem{definition}{Definition}

\numberwithin{equation}{section}

\def\includegraphics{}
\newboolean{publ}

\begin{document}

\centerline{\bf A Novel Subclass of Analytic Functions Specified
by a Family} {\centerline{\bf of Fractional Derivatives in the Complex Domain}
\centerline{}

\centerline{ Zainab Esa$^{1}$, H. M. Srivastava$^{2}$  Adem K{\i}l{\i}cman$^{1}$ and Rabha W. Ibrahim$^{3}$ }

\centerline{$^1$Department of Mathematics, University Putra Malaysia,
43400 UPM Serdang, Selangor, Malaysia}
\centerline{e-mail: esazainab@yahoo.com, \ \ \ akilic@upm.edu.my}
\centerline{}

\centerline{$^2$Department of Mathematics and Statistics,
University of Victoria} 
\centerline{Victoria, British Columbia
V8W 3R4, Canada, \ email: harimsri@math.uvic.ca}
\centerline{}
\centerline{$^3$ Institute of Mathematical Sciences, University of Malaya} 
\centerline{50603 Kuala Lumpur, Malaysia,\ rabhaibrahim@yahoo.com}


\begin{abstract}
\noindent
In this paper, by making use of a certain family of fractional
derivative operators in the complex domain, we introduce and
investigate a new subclass  $\mathcal{P}_{\tau,\mu}(k,\delta,\gamma)$
of analytic and univalent functions in the open unit disk $\mathbb{U}$.
In particular, for functions in the class
$\mathcal{P}_{\tau,\mu}(k,\delta,\gamma)$, we derive sufficient
coefficient inequalities,  distortion theorems involving the
above-mentioned fractional derivative operators, and the radii
of starlikeness and convexity. In addition, some applications of
functions in the class $\mathcal{P}_{\tau,\mu}(k,\delta,\gamma)$
are also pointed out.\\


\end{abstract}

\ifthenelse{\boolean{publ}}{\begin{multicols}{2}}{}

\section{Introduction}

Let  $\mathcal{H}$ be the class of functions which are analytic in the open unit disk
$$\mathbb{U}:= \{z:  z \in \mathbb{C} \quad \text{and} \quad  |z| < 1\}.$$
Also let $\mathcal{H}[a,k]$ denote the subclass of $\mathcal{H}$ consisting of
analytic functions of the form:
$$ f(z)= a + \sum_{j=k}^{\infty}a_{j}z^{j}= a+a_{k}z^{k}+a_{k+1}z^{k+1}+\cdots.$$
We denote by $\mathcal{A}(k)$  the class of functions $f(z)$ {\it normalized} by
\begin{align}\label{+}
f(z)= z + \sum_{\nu=k+1}^{\infty} a_{\nu} z^{\nu} \qquad
(z\in\mathbb{U};\; k \in \mathbb{N}:= \{1,2,3,\cdots\}),
\end{align}
which are analytic in the open unit disk $\mathbb{U}$. In particular, we write
$$\mathcal{A}(1)=: \mathcal{A}.$$
Let $\mathcal{S}(k)$ denote the subclass of
$\mathcal{A}(k)$ consisting of functions which are univalent in $\mathbb{U}$.
Then, by definition, a function $f(z)$ belonging to the univalent function
class $\mathcal{S}(k)$ is said to be a starlike function of order $\alpha\;$
$(0\leqq\alpha <1) $ in $\mathbb{U}$ if and only if
\begin{align}
{\Re} \left(\dfrac{z f^{\prime}(z)}{f(z)}\right) > \alpha \qquad  (z \in \mathbb{U};
\; 0\leqq\alpha <1).
\end{align}
Furthermore, a function $f(z)$ in the univalent function class
$\mathcal{S}(k)$ is said to be a convex function of order
$\alpha\;$ $(0\leqq \alpha < 1)$ in in $\mathbb{U}$ if and only if
\begin{align}\label{Convex}
{\Re} \left(1+ \dfrac{z f^{\prime \prime}(z)}{f^{\prime}(z)}\right) > \alpha
\qquad  (z \in \mathbb{U};
\; 0\leqq\alpha <1).
\end{align}
We denote by $ \mathcal{S}^{*}(k,\alpha)$ and  $\mathcal{K}(k,\alpha)$
the classes of all functions in $ \mathcal{S}(k)$ which are, respectively,
starlike of order $\alpha\; (0 \leqq \alpha<1)$ in
$\mathbb{U}$  and convex  of order $ \alpha\; (0 \leqq \alpha<1)$ in $\mathbb{U}$.

Let $\mathcal{P}(k)$ denote the subclass of $\mathcal{S}(k)$
consisting of functions $f(z)$ which are analytic {\it and} univalent
in $\mathbb{U}$ with negative coefficients, that is, of the form:
\begin{align}\label{-}
 f(z) = z - \sum_{\nu=k+1}^{\infty} a_{\nu} z^{\nu}
\qquad (z\in \mathbb{U};\; a_{\nu} \geqq 0).
\end{align}
For $0\leqq \alpha <1$ and $k \in \mathbb{N}$, we write
\begin{align}
\mathcal{P}^{*}(k,\alpha):= \mathcal{S}^{*}(k,\alpha) \cap  \mathcal{P}(k)
\qquad \text{and}\qquad  \mathcal{L}(k,\alpha)
:= \mathcal{K}(k,\alpha)\cap  \mathcal{P}(k).
\end{align}

Chatterjea \cite{z20} studied the classes $\mathcal{P}^{*}(k,\alpha)$ and
$\mathcal{L}(k,\alpha)$, which are, respectively, starlike and convex of
order $\alpha$ in $\mathbb{U}$. Subsequently, Srivastava {\it et al.} \cite{z21}
observed and remarked that some of the results of
Chatterjea \cite{z20} would follow immediately by trivially setting
$$a_k=0 \qquad (k\in \mathbb{N}\setminus\{1\}=\{2,3,4,\cdots\})$$
in the corresponding earlier results of Silverman \cite[p. 110, Theorem 2; p. 111,
Corollary 2]{z18} (see, for details, \cite[p. 117]{z21}).\\

The {\it modified} convolution of two analytic functions $f(z)$ and $\psi(z)$
in the class $\mathcal{P}(k)$ is defined by \cite{z40}
$$ f * \psi(z):= z- \sum_{\nu=k+1}^{\infty} a_{\nu}\lambda_{\nu}
z^{\nu}=:\psi * f(z),$$
where $f(z)$ is given by \eqref{-} and  $\psi(z)$ is defined as follows:
\begin{equation} \label{psi}
\psi(z)= z- \sum_{\nu=k+1}^{\infty} \lambda_{\nu} z^{\nu}\qquad
(\lambda_{\nu}  \geqq 0;\; k \in \mathbb{N}).
\end{equation}

\begin{definition} \label{def:Int Ser}
{\rm The fractional integral of order $\varsigma$ is defined, for a function $f(z)$ by
\begin{equation}
\mathcal{I}^{  \varsigma}_{z} f(z)
:= \dfrac{1}{\Gamma(\varsigma)} \int_{0}^{z} f(\zeta)
(z- \zeta)^{\varsigma -1} \; {\rm d}\zeta ,
\end{equation}
where $0 \leqq \varsigma < 1,$ the function $f(z)$
is analytic in a simply-connected region
of the complex $z$-plane $\mathbb{C}$ containing the origin
and the multiplicity of $(z-\zeta)^{\varsigma-1}$ is removed by
requiring $\log(z-\zeta)$ to be real when $z-\zeta> 0$.}
\end{definition}

Here, and in what follows, we refer to $\mathcal{I}^{\varsigma}_{z} f(z)$ as
the Srivastava-Owa operator of fractional integral.
Similarly, we have the following definition of
the Srivastava-Owa operator of fractional derivative (see also \cite{KST}).

\begin{definition} \label{def:D sirvastava}
{\rm The fractional derivative of order $\varsigma$ is
defined, for a function $ f(z)$, by
\begin{equation}
\mathfrak{D}^{\varsigma}_{z} f(z)
:= \dfrac{1}{\Gamma(1-\varsigma)}\;\dfrac{\rm d}{{\rm d}z}\left\{
\int_{0}^{z} f(\zeta) (z- \zeta)^{-\varsigma} {\rm d}\zeta \right\},
\end{equation}
where $0 \leqq \varsigma < 1,$ the function $f(z)$ is analytic in a
simply-connected region of the complex $z$-plane $\mathbb{C}$ containing the origin
and the multiplicity of $(z-\zeta)^{-\varsigma}$ is removed as in Definition \ref{def:Int Ser}}.
\end{definition}

Now, by using Definition \ref{def:D sirvastava}, the Srivastava-Owa fractional
derivative of order  $n+\varsigma$ can easily be defined as follows:
\begin{align}\label{2}
\mathfrak{D}^{n+\varsigma}_{z} f(z) := \dfrac{{\rm d}^{n}}{{\rm d}z^{n}}
\left\{\mathfrak{D}^{\varsigma}_{z} f(z)\right\}
\qquad
(0 \leqq \varsigma < 1;\;  n \in \mathbb{N}_{0}
:=\mathbb{N}\cup\{0\}=\{0,1,2,3,\cdots\}),
\end{align}
which readily yields

$$\mathfrak{D}^{0+\varsigma}_{z} f(z)
=\mathfrak{D}^{\varsigma}_{z} f(z) \qquad  \text{and} \qquad
\mathfrak{D}^{1+\varsigma }_{z} f(z) = \dfrac{\rm d}{{\rm d}z}
\left\{\mathfrak{D}^{\varsigma}_{z} f(z)\right\} \qquad ( 0 \leqq \varsigma < 1).
$$

Recently, by applying the Srivastava-Owa definition (\ref{2}), Tremblay \cite{z30}
introduced and studied an interesting fractional derivative operator $\; {\mathfrak T}^{\tau,\mu},$
which was defined in the complex domain and whose properties in several spaces were discussed systematically
(see, for details, \cite{2} and \cite{z30}).

\begin{definition}\label{T3}
{\rm For $0 < \tau \leqq 1,\,$  $0< \mu  \leqq 1\,$ and $\, 0 \leqq \tau-\mu  < 1,\;$
the Tremblay operator $\; {\mathfrak T}^{\tau,\mu}\;$  of a function
$\; f \in \mathcal{A}\;$ is defined for all  $ \, z \in \mathbb{U}\;$ by}
\begin{equation}
{\mathfrak T}^{\tau,\mu}\;f(z) :=  \frac{\Gamma(\mu)}{\Gamma(\tau )} \; z^{1-\mu}\;
\mathfrak{D}^{\tau-\mu}_{z}\; z^{\tau-1} \;f(z) \qquad (z \in \mathbb{U}) .\label{de:modify T}
\end{equation}
\end{definition}

In the special case when $\tau= \mu=1$ in \eqref{de:modify T}, we have
\begin{align}
{\mathfrak T}^{1,1}f(z) = f(z).
\end{align}
We note also that $\;\mathfrak{D}^{\tau-\mu}_{z}\;$ represents a
Srivastava-Owa operator of fractional derivative of order
$\tau-\mu$ $\; (0 \leqq \tau-\mu <1),\;$ which is given by Definition \ref{def:D sirvastava}.\\

The main purpose of this paper is to present coefficient inequalities and
coefficient estimates, distortion theorems, and the radii of  starlikeness
and convexity, for functions belonging to the class $\mathcal{P}_{\tau,\mu}(k,\delta,\gamma)$
which we introduce in Section 2 below.
We also consider some other interesting results involving closure and
convolution of functions in the class $\mathcal{P}_{\tau,\mu}(k,\delta,\gamma)$.

\section{A Set of Main Results}

In this section, we define a new analytic class $\mathcal{P}_{\tau,\mu}(k,\delta,\gamma)$
by considering the fractional derivative operator given by Definition \ref{T3}
and establish a sufficient condition for a function $f(z) \in \mathcal{P}(k)$ to be in
the function class $\mathcal{P}_{\tau,\mu}(k,\delta,\gamma)$. The following two lemmas
will be needed in our investigation.

\begin{lemma}\label{lemma1}
Let the function $f(z)$ defined by $\eqref{-}$ belong to the class
$\mathcal{P}(k)\;\;( k\in \mathbb{N})$. Then
$${\mathfrak T}^{\tau,\mu}\; f(z) =\frac{\tau}{\mu} z- \sum_{\nu=m+1}^{\infty}
\frac{\Gamma(\nu+\tau)\Gamma(\mu)}{\Gamma(\nu++\mu)\Gamma(\tau)}\;  a_{\nu} z^{\nu },$$
where $\; 0 < \tau \leqq 1,\;$ $\; 0 < \mu \leqq 1\;$ and $\; 0 \leqq \tau- \mu < 1$.
\end{lemma}

\noindent
\textbf{Proof.} By using Definition \ref{T3} and Definition \ref{def:D sirvastava},
we find for $\;z\in \mathbb{U}\;$ that
\begin{align*}
{\mathfrak T}^{\tau,\mu} f(z) &= \frac{\Gamma(\mu)}{\Gamma(\tau)}
\; z^{1-\mu} \;\mathfrak{D}^{\tau-\mu}_{z} z^{\tau-1}\; f(z) \\
&=\frac{\Gamma(\mu)}{\Gamma(\tau)} \; z^{1-\mu} \;\mathfrak{D}^{\tau-\mu}_{z} \;z^{\tau-1}
\left(z- \sum_{\nu=m+1}^{\infty} a_{\nu} z^{\nu}\right)\\
&= \frac{\Gamma(\mu)}{\Gamma(\tau)} \;z^{1-\mu}\; \mathfrak{D}^{\tau-\mu}_{z}
\left(z^{\tau}- \sum_{\nu=k+1}^{\infty} a_{\nu} z^{\nu+\tau-1}\right)\\
&=\frac{\Gamma(\mu)}{\Gamma(\tau)}\; z^{1-\mu} \;
\left(\frac{\Gamma(\tau+1)}{\Gamma(\mu+1)}\;
z^{\mu}- \sum_{\nu=k+1}^{\infty}\frac{\Gamma(\nu+\tau)}
{\Gamma(\nu+\mu)} \; a_{\nu} z^{\nu+\mu-1}\right)\\
&=\frac{\tau}{\mu}\; z- \sum_{\nu=k+1}^{\infty}
\frac{\Gamma(\nu+\tau)\Gamma(\mu)}
{\Gamma(\nu+\mu)\Gamma(\tau)}\;  a_{\nu} z^{\nu },
\end{align*}
which proves Lemma \ref{lemma1}.

\begin{lemma}\label{lemma2}
Let the function $f(z)$ defined by $\eqref{-}$
belong to class $\mathcal{P}(k)\;\;( k\in \mathbb{N})$. Then
\begin{align*}
\bigl({\mathfrak T}^{\tau,\mu}\; f(z)\bigr)^{\prime}
= \frac{\tau}{\mu} - \sum_{\nu=k+1}^{\infty}
\frac{\nu\Gamma(\nu+\tau)\Gamma(\mu)}
{\Gamma(\nu+\mu)\Gamma(\tau)}\;  a_{\nu} z^{\nu-1 }
\qquad (z\in \mathbb{U}),
\end{align*}
where $\; 0 < \tau \leqq 1,\;$ $\; 0 < \mu \leqq 1\;$
and $\;0 \leqq \tau- \mu < 1$.
\end{lemma}

\noindent
\textbf{Proof.}
By using Lemma \ref{lemma1} and the definition \ref{2}, we have
\begin{align*}
\dfrac{\rm d}{{\rm d}z} \bigl\{{\mathfrak T}^{\tau,\mu}\; f(z)\bigr\}
&= \dfrac{\rm d}{{\rm d}z}\left\{\frac{\Gamma(\mu)}{\Gamma(\tau)}\;
z^{1-\mu} \; \mathfrak{D}^{\tau-\mu}_{z} \;z^{\tau-1} \;f(z)\right\}\\
&=\frac{\tau}{\mu} - \sum_{\nu=k+1}^{\infty} \frac{\nu \Gamma(\nu+\tau)\Gamma(\mu)}
{\Gamma(\nu+\mu)\Gamma(\tau\alpha)} \; a_{\nu} z^{\nu-1 } \qquad (z\in \mathbb{U}),
\end{align*}
which evidently completes the proof of Lemma \ref{lemma2}

By employing Lemma \ref{lemma1} and  Lemma \ref{lemma2},
we now introduce a new class
$\mathcal{P}_{\tau,\mu}(k,\delta,\gamma)$ of analytic
functions in $\mathbb{U}$  as follows.

\begin{definition}\label{def4}
{\rm Let $\; 0 <  \tau  \leqq 1,\;$ $\;0 <  \mu  \leqq 1,\;$
$\;0 \leqq \delta < 1\;$ and $\;0 \leqq \gamma < 1$.
A function $ f(z)$  belonging to the analytic function class
$\mathcal{P}(k)$ is said to be in the
class $ \mathcal{P}_{\tau,\mu}(k,\delta,\gamma)$ if and only if
\begin{align}\label{Inq:1}
{\Re}\left(\dfrac{\Gamma(\mu+1)\Gamma(\tau)}{\Gamma(\tau+1)\Gamma(\mu)}\;
z^{-1} \left[(1-\delta)\, {\mathfrak T}^{\tau,\mu}\; f(z)
+ z\delta\bigl( {\mathfrak T}^{\tau, \mu}\;
f(z)\bigr)^{\prime}\right]\right) > \gamma  \qquad
(z \in \mathbb{U};\; \tau-\mu + \gamma <1),
\end{align}
where ${\mathfrak T}^{\tau, \mu}$ is the fractional derivative operator in
the complex domain in Definition $\ref{T3}$.}
\end{definition}

\subsection{A Theorem on Coefficient Bounds}

\begin{theorem}\label{Th:1}
Let the function  $f(z)$ be given by $\eqref{-}$. Then $f(z)$ belongs
to the class $\;\mathcal{P}_{\tau,\mu}(k,\delta,\gamma)\;$
if and only if
\begin{align}\label{Inq10}
\sum_{\nu =k+1 }^{ \infty} \dfrac{(1+\delta \nu - \delta ) \,
\Gamma(\nu+\tau)\Gamma(\mu+1)}{\Gamma(\nu+\mu)\Gamma(\tau+1)} \, a_{\nu}
\leqq 1-\gamma \qquad (a_{\nu}\geqq 0; \, 0 \leqq \gamma < 1).
\end{align}
The result $\eqref{Inq10}$ is sharp and the extremal function $f(z)$ is given by by
\begin{align}\label{sharp}
f(z) = z- \dfrac{(1- \gamma) (\mu+1)_{k}}{(1+\delta k) (\tau+1)_{k}}\;z^{k+1}
\qquad (k \in \mathbb{N}),
\end{align}
where $\;(\lambda)_{k}\;$ denotes the Pochhammer symbol defined by

\begin{equation*}
(\lambda)_{k}:=\dfrac{\Gamma(\lambda+k)}{\Gamma(\eta)}=
\begin{cases}
1&  \qquad (k = 0)\\
\\
\lambda(\lambda+1) \cdots (\lambda+k-1) &  \qquad (k \in  \mathbb{N}).
\end{cases}
\end{equation*}
\end{theorem}

\noindent
\textbf{Proof.}
Supposing first that $\;f(z)\in\mathcal{P}_{\tau,\mu}(k,\delta,\gamma),\; $ we find
from Definition \ref{T3} in conjunction with with Lemmas \ref{lemma1} and  \ref{lemma2} that
\begin{align}\label{Inq:1}
\Re \left(1- \sum_{\nu =k+1 }^{ \infty} \dfrac{(1+\delta \nu- \delta ) \,
\Gamma(\nu+\tau)\Gamma(\mu+1)}{\Gamma(\nu+\mu)\Gamma(\tau+1)}
\,a_{\nu} \, z^{\nu-1}  \right)  > \gamma.
\end{align}
If we choose $z$ to be real and let $z \rightarrow -1$, we have
\begin{align*}
1-\sum_{\nu =k+1 }^{ \infty} \dfrac{(1+\delta \nu- \delta )  \,
\Gamma(\nu+\tau)\Gamma(\mu+1)}
{\Gamma(\nu+\mu)\Gamma(\tau+1)}\; a_{\nu}
\geqq \gamma \qquad (0 < \tau \leqq 1;\; 0 < \mu \leqq 1),
\end{align*}
which readily yields the inequality \eqref{Inq10} of Theorem \ref{Th:1}.

Conversely, by assuming that the inequality \eqref{Inq10}
is true, we let $|z|=1$. We then obtain

\begin{align}\label{Inq:1}
&\left|\dfrac{\Gamma(\mu+1)\Gamma(\tau)}{\Gamma(\tau+1)\Gamma(\mu)}\;z^{-1}
\left[  (1-\delta) {\mathfrak T}^{\tau-\mu}\; f(z) + \delta z
\left({\mathfrak T}^{\tau-\mu} \;f(z)\right)^{\prime}  \right] -1 \right|\notag \\
&\qquad \quad = \left|- \sum_{\nu =k+1 }^{ \infty} \dfrac{(1+\delta \nu- \delta )  \,
\Gamma(\nu+\tau)\Gamma{(\mu+1)}}{\Gamma(\nu+\mu)\Gamma(\tau+1)}
\;a_{\nu}\;  z^{\nu-1}\right|\notag \\
&\qquad \quad \leqq \sum_{\nu =k+1 }^{ \infty} \dfrac{(1+\delta \nu- \delta )  \,
\Gamma(\nu+\tau)\Gamma{(\mu+1)}}{\Gamma(\nu+\mu)
\Gamma(\tau+1)}\; a_{\nu}\;|z|^{\nu-1}\notag \\
&\qquad \quad \leqq  1- \gamma,
\end{align}
which shows that the function $f(z) $ is in the class
$\;\mathcal{P}_{\tau,\mu}(k,\delta,\gamma)$.

Finally, it is easily verified that the result is sharp
for the function $f(z)$ given by \eqref{sharp}.

\begin{corollary}
Let the function $\;f(z)\;$ given by $\eqref{-}$ be in the class
$\; \mathcal{P}_{\tau,\mu}(k,\delta,\gamma)$. Then
\begin{align}\label{Exm1}
 a_{k+1} \leqq \dfrac{(1- \gamma) (\mu+1)_{k} }{(1+\delta k) (\tau+1)_{k}}
\end{align}
$$(k= \mathbb{N}\setminus\{1\}=\{2,\ 3,\ 4,\cdots\};\; 0< \mu  \leqq 1;\;
0  < \tau  \leqq  1;\; 0 \leqq  \gamma <  1;\; 0 \leqq \delta < 1).$$
\end{corollary}

\begin{corollary}\label{Corollary1}
The function $\; f(z) \in \mathcal{P}(k)\;$ is in the class
$\mathcal{P}_{1,1}(k,\delta,\gamma)$ if and only if
\begin{align}
\sum_{\nu=k+1}^{\infty} \left( 1+ \delta\nu - \delta \right)
a_{\nu} \leqq 1-\gamma
\qquad ( 0 \leqq \gamma < 1;\; 0 \leqq \delta \leqq 1).
\end{align}
\end{corollary}

Corollary \ref{Corollary1} was given  by Altinta\c{s} {\it et al.}
\cite{z41}. In particular, it was given earlier for $k=1$
by Bhoosnurmath and Swamy \cite{z42} for $k=1$
and by Silverman \cite{z18} for $k=\delta =1$.

\subsection{Distortion Theorems}
\begin{theorem}\label{TH1}
Let the function $\;f(z)\;$ belong to the class  $\;\mathcal{P}_{\tau,\mu}(k,\delta,\gamma)$. Then
\begin{align}\label{Dis1}
&\dfrac{\beta}{\alpha}\; |z|\left(1 - |z|^{k}
\dfrac{(1-\gamma)(\beta+1)_{k}(\mu+1)_{k}}{(1+\delta k)\, (\alpha +1)_{k}(\tau +1)_{k} }\right)\notag \\
&\qquad \quad \leqq \bigl| {\mathfrak T}^{\beta,\alpha}\;f(z)\bigr| \notag \\
&\qquad \quad \leqq   \dfrac{\beta}{\alpha} \;|z|\left(1+|z|^{k}\dfrac{(1-\gamma)(\beta+1)_{k}(\mu+1)_{k}}
{(1+\delta k)\, (\alpha+1)_{k}(\tau+1)_{k}}\right)
\end{align}
$$(z \in \mathbb{U};\; 0 < \mu \leqq 1;\; 0 < \beta \leqq 1;\; k \in \mathbb{N}_{0}.$$
\end{theorem}

\noindent
\textbf{Proof.}
By hypothesis, the function $\;f(z)\;$ belongs to the class
$\mathcal{P}_{\tau,\mu}(k,\delta,\gamma)$. Thus, clearly, we find from the inequality
\eqref{Inq10} in Theorem \ref{Th:1} that
\begin{align}
\dfrac{(1+\delta k) \, (\tau+1)_{k}}{(\mu+1)_{k}}
\sum_{\nu=k+1}^{\infty}  a_{\nu} \leqq  \sum_{\nu =k+1}^{\infty}
\dfrac{(1+\delta\nu- \delta)) \, \Gamma(\nu+\tau)\Gamma{(\mu+1)}}{\Gamma(\nu+\mu)\Gamma(\tau+1)}\; a_{\nu},
 \end{align}
which leads us to
\begin{align}\label{a}
\sum_{\nu=k+1}^{\infty} a_{\nu} \leqq  \dfrac{(1-\gamma)(\mu +1)_{k}}{(1+\delta k ) \, (\tau +1)_{k}}
\end{align}
$$(0<\mu\leqq 1;\; 0<\tau\leqq 1;\;  0\leqq \delta<1;\;0\leqq \gamma<1;\; k\in \mathbb{N}).$$

Next, by the definition \eqref{de:modify T} and from \eqref{a}, we have
\begin{align*}
{\mathfrak T}^{\beta,\alpha} \;f(z)
&= \dfrac{\Gamma(\alpha)}{\Gamma(\beta)} \;z^{1-\alpha}
D^{\beta-\alpha}_{z}\; z^{\beta-1}\; f(z)  \\
&=\dfrac{\beta}{\alpha} \left(z - \sum_{\nu=k+1}^{\infty} \dfrac{\Gamma(\nu+\beta)\Gamma(\alpha+1)}
{\Gamma(\nu+\alpha) \Gamma(\beta+1)} \;a_{\nu} \; z^{\nu}\right)\\
&= \dfrac{\beta}{\alpha} \left(z- \sum_{\nu=k+1}^{\infty} \omega(\nu) a_{\nu} \; z^{\nu} \right),
\end{align*}
where
\begin{align*}
\omega(\nu)=\dfrac{(\beta +1)_{\nu-1} }{(\alpha+1)_{\nu-1}}
\qquad ( \nu= k+1,k+2,k+3,\cdots).
\end{align*}
Since the function $\omega(\nu)$ can be seen to be non-increasing, we get
\begin{align}\label{C}
0< \omega(\nu) \leqq  \omega(k+1)
= \dfrac{(\beta+1)_{k}}{(\alpha+1)_{k}}.
\end{align}
Thus, from the inequalities \eqref{C} and \eqref{a}, we find that
\begin{align*}
\bigl| {\mathfrak T}^{\beta,\alpha}\;f(z)\bigr|
&\geqq \dfrac{\beta}{\alpha} \left(|z|  -
\left| \omega(k+1) \sum_{\nu=k+1}^{\infty} a_{\nu} z^{\nu}\right|\right)\\
& \geqq \dfrac{\beta}{\alpha} \left(|z|  - |z|^{\nu} \omega(k+1)
\sum_{\nu=k+1}^{\infty} a_{\nu}\right)\\
&\geqq \dfrac{\beta}{\alpha} |z| \left(1 - |z|^{k}
\dfrac{(\beta+1)_{k}\,(1-\gamma)(\mu+1)_{k}}
{(\alpha+1)_{k}(1+\delta k)\, (\tau+1)_{k}}\right),
\end{align*}
which proves the first part of the inequality \eqref{Dis1}. In a similar manner,
we can prove the second part of the inequality \eqref{Dis1}.

\begin{theorem}\label{Th5}
Let the function $\;f(z)\;$ be in the class
$\;\mathcal{P}_{\tau,\mu}(k,\delta,\gamma)$. Then
\begin{align}\label{fDis1}
|z| -|z|^{k+1}\; \dfrac{(1- \gamma) (\mu +1)_{k}}{(1 +\delta k) (\tau+1)_{k}}
\leqq  \,| f(z) | \leqq   |z| +  |z|^{k+1}\;
\dfrac{(1- \gamma) (\mu +1)_{k}}{(1 +\delta k) (\tau+1)_{k}}.
\end{align}
\end{theorem}

\noindent
\textbf{Proof.}
By using  the same method as in Theorem \ref{TH1}, we have
\begin{align}\label{f1}
|f(z)|&\leqq  |z| + \sum_{\nu=k+1}  a_{\nu} |z|^{\nu}  \nonumber \\
& \leqq |z| + |z|^{k+1}\;  \sum_{\nu=k+1}  a_{\nu} \nonumber \\
& \leqq |z| +  |z|^{k+1}\;\dfrac{(1- \gamma) \Gamma(k+1+\mu)\Gamma(\tau+1)}
{(1+\delta k) \Gamma(k+1+\tau) \Gamma(\mu+1)} \nonumber \\
&= |z| + |z|^{k+1}\;\dfrac{(1- \gamma)  (\mu+1)_{k}}{(1+\delta k) \,(\tau+1)_{k}}
\end{align}
and

\begin{align}\label{f2}
|f(z)| &\geqq |z| - |z|^{k+1}\;\sum_{\nu=k=1}  a_{\nu} \nonumber \\
&\geqq  |z| - |z|^{k+1}\l\dfrac{(1- \gamma)  (\mu+1)_{k}}{(1+\delta k) \,(\tau+1)_{k}}.
\end{align}
Consequently, from \eqref{f1} and \eqref{f2},  we immediately get the
inequality \eqref{fDis1} of Theorem \ref{Th5}.\\

Upon setting $\tau=\mu=1$ in Theorem \ref{Th5}, we obtain the following corollary.

\begin{corollary}
If $\; f(z) \in \mathcal{P}_{0}(k,\delta,\gamma)=: \mathcal{P}(k,\delta,\gamma),\;$ then
\end{corollary}
\begin{align*}
|z| - |z|^{k+1} \; \left(\dfrac{1-\gamma}{1+\delta k} \right) \leqq
\left|  f(z)\right| \leqq   |z| + |z|^{k+1} \;
\left(\dfrac{1-\gamma}{1+\delta k}\right)
\end{align*}
$$(z\in \mathbb{U};\; 0 \leqq \gamma< 1;\; 0 \leqq \delta < 1;\;  k \in \mathbb{N}).$$

Moreover, if $\;\tau=\mu=1\;$ and $\;k=1\;$ in Theorem \ref{Th5}, then we have
the following known result (see \cite{z42}).

\begin{corollary}
If $\;f(z) \in \mathcal{P}_{0}(1,\delta,\gamma)=: \mathcal{P}(\delta,\gamma),\;$ then
\begin{align*}
|z| - |z|^{2}\;\left(\dfrac{1-\gamma}{1+\delta}\right) \leqq
\left|f(z)\right| \leqq   |z| + |z|^{2} \; \left(\dfrac{1-\gamma}{1+\delta}\right)
\qquad (z\in \mathbb{U}).
\end{align*}
\end{corollary}

\section{Radii of Starlikeness and Convexity}

\begin{theorem}\label{Theorem4}
If the function $\;f(z) \in  \mathcal{P}_{\tau,\mu}(k,\delta,\gamma),\;$ then
$\;f(z) \in  \mathcal{P}_{\tau,\mu}^{*}(k,\delta,\gamma)\;$ in the disk $\;|z| < r_{1},\;$
where
$$r_{1}:= \inf_{\nu  \geqq k+1} \left\{
\dfrac{(1-\alpha)(1+\delta\nu- \delta)\Gamma(\nu+\tau)\Gamma(\mu+1)}
{(\nu-\alpha)(1-\gamma) \Gamma(\nu+\mu)\Gamma(\tau+1)}\right\} ^{1/(\nu-1)}.$$
\end{theorem}

\noindent
\textbf{Proof.}
We must show that the condition in \eqref{Convex} holds true. Indeed, since
\begin{align}
\left| \dfrac{z f^{\prime }(z)}{f(z)} - 1 \right|
\leqq  \frac{ z- \sum\limits_{\nu=k+1}^{\infty} \nu a_{\nu} |z|^{\nu}}
{z-\sum\limits_{\nu=2}^{\infty}  a_{\nu} |z|^{\nu}}  \leqq 1- \alpha
\qquad (z\in \mathbb{U})
\end{align}
and
\begin{align}
{\sum_{\nu=k+1}^{\infty} (\nu - \alpha) a_{\nu} |z|^{\nu-1}} \leqq   1-\alpha
\qquad (z\in \mathbb{U}),
\end{align}
we find that
$$\frac{(\nu-\alpha) |z|^{\nu-1}}{1-\alpha } \leqq  \dfrac{(1+\delta\nu-
\delta)\Gamma(\nu+\tau)\Gamma(\mu+1)}{(1-\gamma)
\Gamma(\nu+\mu)\Gamma(\tau+1)}\qquad (\nu\geqq k+1),$$
that is, that
$$|z| \leqq \left\{ \dfrac{(1-\alpha)
(1+\delta\nu- \delta)\Gamma(\nu+\tau)\Gamma(\mu+1)}
{(\nu-\alpha)(1-\gamma) \Gamma(\nu+\mu)\Gamma(\tau+1)}\right\}^{1/(\nu-1)},$$
which proves Theorem \ref{Theorem4}.

\begin{corollary}\label{Coro;S}
If the function  $\;f(z) \in  \mathcal{P}_{1,1}(k,\delta,\gamma),\;$ then $f(z)$ is
starlike of order $\alpha$ in the disk $\; |z| < r_{2},\;$
where
$$ r_{2}:= \inf_{\nu  \geqq k} \left\{\dfrac{(1-\alpha) (1+\delta\nu- \delta)}
{(\nu-\alpha)(1-\gamma)}\right\} ^{1/(\nu-1)}.$$
\end{corollary}

In its special case when $k=1$, Corollary \ref{Coro;S} was proven
by Altinta\c{s} {\it et al.} \cite{z41}. Moreover,
Corollary \ref{Coro;S} was given earlier
by Bhoosnurmath and Swamy \cite{z42} for $k=1$ and $\alpha =0$,
and by Silverman  \cite{z18}
when $k = 1$ and $\delta=\gamma = 0$.

\begin{theorem}\label{Theorem5}
If the function  $\;f(z) \in  \mathcal{P}_{\tau,\mu}(k,\delta,\gamma),\;$ then
$ f(z) \in \mathcal{K}_{\tau,\mu}(k,\delta,\gamma)$ in the disk $\;|z| < r_{3},\;$
where
$$ r_{3}:= \inf_{\nu  \geqq k+1} \left\{\dfrac{(1-\alpha)(1+\delta\nu- \delta)
\Gamma(\nu+\tau)\Gamma(\mu+1)}{\nu(\nu-\alpha)(1-\gamma)\Gamma(\nu+\mu)
\Gamma(\tau+1)}\right\}^{1/(\nu-1)}.$$
\end{theorem}

\noindent
\textbf{Proof.}
For the function $ f(z) $ given by \eqref{-}, we must show that
$$\left| \dfrac{zf^{\prime \prime}(z)}{f^{\prime}(z)}\right|
\leqq 1- \alpha \qquad (z \in \mathbb{U}).$$
First of all, we find from  \eqref{-} that
\begin{align*}
\left| \dfrac{zf^{\prime \prime}(z)}{f^{\prime}(z)} \right|
&= \left| \dfrac{-\sum\limits_{\nu=k+1}^{\infty} \nu(\nu-1) a_{\nu}
z^{\nu-1}}{1- \sum\limits_{\nu=k+1}^{\infty} \nu a_{\nu} z^{\nu-1}} \right|\notag \\
&\leqq \dfrac{\sum\limits_{\nu=k+1}^{\infty} \nu(\nu-1) a_{\nu}
|z|^{\nu-1}}{1-\sum\limits_{\nu=k+1}^{\infty} \nu a_{\nu} |z|^{\nu-1}}\notag \\
& \leqq 1- \alpha \qquad (z \in \mathbb{U}),
\end{align*}
if
\begin{align}
\sum_{\nu=k+1}^{\infty} \nu(\nu-1) a_{\nu} |z|^{\nu-1}
\leqq (1-\alpha) \left(1-\sum_{\nu=k+1}^{\infty} \nu a_{\nu} |z|^{\nu-1}\right)
\qquad (z \in \mathbb{U}),
\end{align}
that is, if
\begin{align}\label{L}
\sum_{\nu=k+1}^{\infty} \nu(\nu-\alpha) a_{\nu} |z|^{\nu-1} \leqq 1-\alpha
\qquad (z \in \mathbb{U}).
\end{align}

From the last inequality \eqref{L}, together with Theorem \ref{Th:1}, we thus
find that
\begin{align*}
\frac{\nu(\nu-\alpha) |z|^{\nu-1}}{(1-\alpha)} \leqq  \dfrac{(1+\delta \nu -\delta)
\Gamma(\nu+\tau)\Gamma(\mu+1)}{(1-\gamma)\Gamma(\nu+\mu)\Gamma(\tau+1)}
\qquad (\nu \geqq k+1),
\end{align*}
that is, that

\begin{align*}
|z| \leqq
 \left\{\dfrac{(1-\alpha)(1+\delta \nu -\delta)
\Gamma(\nu+\tau)\Gamma(\mu+1)}{\nu(\nu-\alpha)
(1-\gamma)\Gamma(\nu+\mu)\Gamma(\tau+1)}\right\}^{1/(\nu-1)},
\end{align*}
which evidently proves Theorem \ref{Theorem5}.

\begin{corollary}\label{Coro;Conv}
If the function  $\;f(z) \in  \mathcal{P}_{0}(k,\delta,\gamma),\;$ then
$f(z)$ is convex of order $\alpha$ in  the disk $\;|z| < r_{4},\;$
where
$$ r_{4}:= \inf_{\nu \geqq k+1} \left\{\dfrac{(1-\alpha)
(1+\delta \nu - \delta)}{\nu(\nu-\alpha)(1-\gamma)}\right\}^{1/(\nu-1)}.$$
\end{corollary}

For $ k=1$, Corollary \ref{Coro;Conv} was proved by
Altinta\c{s} {\it et al.} \cite{z41}.
Further special cases of Corollary \ref{Coro;Conv}
were given earlier by  Bhoosnurmath and Swamy  \cite{z42}
when $k=1$ and $\alpha=0$, and by Silverman \cite{z18} for
$ k=1$ and $ \delta= \gamma=0$.

\section{Further Results for the Function Class
$\;\mathcal{P}_{\tau,\mu}(k,\delta,\gamma)$}

In this section, we prove some results for the
closure of functions and the convolution of
functions in the class
$\;\mathcal{P}_{\tau,\mu}(k,\delta,\gamma)$.

\begin{theorem}\label{Th.a1,a2}
Let each of the functions $f_{1}(z)$ and $f_{2}(z)$ given by
$$f_{1}(z)= z- \sum_{\nu=k+1}^{\infty} a_{\nu,1} z^{\nu}
\qquad (a_{\nu,1}\geqq 0;\; k \in \mathbb{N})$$
and
$$f_{2}(z)= z- \sum_{\nu=k+1}^{\infty} a_{\nu,2} z^{\nu}
\qquad (a_{\nu,2}\geqq 0;\; k \in \mathbb{N})$$
be in the class $\;\mathcal{P}_{\tau,\mu}(k,\delta,\gamma)$.
Then the function $\Phi(z)$ given by
$$\Phi(z) = z- \frac{1}{2} \sum_{\nu=k+1}^{\infty} (a_{\nu,1}
+ a_{\nu,2} )z^{\nu}$$
is also in the class $\;\mathcal{P}_{\tau,\mu}(k,\delta,\gamma)$.
\end{theorem}

\noindent
\textbf{Proof.}
By the hypothesis that each of the functions $f_{1}(z)$ and
$f_{2}(z)$ is in class $ \mathcal{P}_{\tau,\mu}(k,\delta,\gamma)$, we get

\begin{align}\label{Inq:1}
\sum_{\nu =k+1 }^{ \infty} \dfrac{(1+\delta \nu - \delta) \,
\Gamma(\nu+\tau)\Gamma{(\mu+1)}}{\Gamma(\nu+\mu)\Gamma(\tau+1)}\;
a_{\nu,1} \leqq   1- \gamma
\end{align}
and
\begin{align}
\sum_{\nu =k+1 }^{ \infty} \dfrac{(1+\delta \nu - \delta) \,
\Gamma(\nu+\tau)\Gamma{(\mu+1)}}{\Gamma(\nu+\mu)\Gamma(\tau+1)} \;
a_{\nu,2} \leqq   1 - \gamma ,
\end{align}
so that, obviously,
\begin{align}
\frac{1}{2} \sum_{\nu =k+1 }^{\infty}& \dfrac{(1+\delta \nu -\delta) \,
\Gamma(\nu+\tau)\Gamma{(\mu+1)}}{\Gamma(\nu+\mu)\Gamma(\tau+1)} \;
( a_{\nu ,1} + a_{\nu,2} )
\leqq 1- \gamma,
\end{align}
which proves the assertion that
$\Phi(z) \in  \mathcal{P}_{\tau,\mu}(k,\delta,\gamma).$

\begin{theorem}\label{Theorem7}
Let the functions $f_{j}(z)$ $\;(j=1,\cdots,p)$ defined by
$$f_{j}(z)= z- \sum_{\nu=k+1}^{\infty} a_{\nu,j} z^{\nu}
\qquad (a_{\nu,j} \geqq 0;\; k \in \mathbb{N})$$
be in the class $\;\mathcal{P}_{\tau,\mu}(k,\delta,\gamma)$.
Then the function $\Theta(z)$ defined by
\begin{equation}\label{Theta}
\Theta(z):= \sum_{j=1}^{p} q_{j} f_{j}(z) \qquad (q_{j} \geqq 0)
\end{equation}
is also  in the class $\;\mathcal{P}_{\tau,\mu}(k,\delta,\gamma),\;$ where
$$ \sum_{j=1}^{\infty} q_{j} = 1\qquad (q_{j} \geqq 0).$$
\end{theorem}

\noindent
\textbf{Proof.}
By the definition \eqref{Theta} of the function $\Theta(z)$, we have
\begin{align*}
\Theta(z) &=\sum_{j=1}^{p} q_{j}\left( z-
\sum_{\nu=k+1}^{\infty} a_{\nu,j} z^{\nu}\right)\\
&=\sum_{j=1}^{p} q_{j} z -\sum_{\nu=k+1}^{\infty}
\left( \sum_{j=1}^{p} q_{j} a_{\nu,j} z^{\nu}\right)\\
&=z -\sum_{\nu=k+1}^{\infty} \left( \sum_{j=1}^{p} q_{j}
a_{\nu,j} \right)z^{\nu}.
\end{align*}
Since $f_{j}(z)  \in  \mathcal{P}_{\tau,\mu}(k,\delta,\gamma)$
$\;(j=1,\cdots,p)$, we also have
\begin{align*}
\sum_{\nu =k+1 }^{ \infty} \dfrac{(1+\delta \nu - \delta) ) \,
\Gamma(\nu+\tau)\Gamma{(\mu+1)}}{\Gamma(\nu+\mu)\Gamma(\tau+1)} a_{\nu,j}
\leqq   1- \gamma \qquad (j=1,\cdots,p).
\end{align*}
The remainder of the proof of Theorem \ref{Theorem7} (which is based essentially
upon Theorem \ref{Th:1}) is fairly straightforward
and is, therefore, omitted here.

\begin{theorem}\label{Theorem8}
Let the function $f(z)$ given by $\eqref{-}$ and the function $\hslash(z)$ defined by
$$\hslash(z)= z- \sum_{\nu=k+1}^{\infty} \lambda_{\nu} z^{\nu}
\qquad (\lambda_{\nu}\geqq 0;\; k \in \mathbb{N})$$
be in the same class  $\;\mathcal{P}_{\tau,\mu}(k,\delta,\gamma)$.
Then the function
$\Delta(z)$ defined by
\begin{align*}
\Delta(z) &= (1-\eta) f(z)+ \eta \hslash(z)\\
&= z  - \sum_{\nu=k+1}^{\infty} \rho_{\nu} z^{\nu}.
\end{align*}
is also in the class  $\;\mathcal{P}_{\tau,\mu}(k,\delta,\gamma)$.
\end{theorem}

\noindent
\textbf{Proof.}
In view of the hypotheses of Theorem {Theorem8}, we find by using Theorem \ref{Th:1} that
\begin{align*}
\sum_{\nu =k+1}^{\infty} \dfrac{(1+\delta \nu - \delta) \,
\Gamma(\nu+\tau)\Gamma(\mu)}{\Gamma(\nu+\mu)\Gamma(\tau)} \;\rho_{\nu}
&= (1-\eta)\sum_{\nu =k+1}^{\infty} \dfrac{(1+\delta \nu - \delta) \,
\Gamma(\nu+\tau)\Gamma(\mu)}{\Gamma(\nu+\mu)\Gamma(\tau)} \;a_{\nu} z^{\nu-1}\\
&\qquad \quad +
\eta \sum_{\nu =k+1}^{\infty}
\dfrac{(1+\delta \nu - \delta) \,
\Gamma(\nu+\tau)\Gamma(\mu)}{\Gamma(\nu+\mu)\Gamma(\tau)}\; \lambda_{\nu} z^{\nu-1}\\
& \leqq (1-\eta)(1-\gamma)+ \eta (1-\gamma)\\
& \leqq 1-\gamma.
\end{align*}
Hence $\;\Delta(z) \in  \mathcal{P}_{\tau,\mu}(k,\delta,\gamma)$.

\begin{theorem}\label{Theorem9}
Let the function $f(z)$ given by $\eqref{-}$ and the function $\psi(z)$ defined by
$\eqref{psi}$ be in the class $\;\mathcal{P}_{\tau,\mu}(k,\delta,\gamma)$.
Then the function $\Omega(z)$ given by the following modified Hadamard product$:$
$$\Omega(z):=f* \psi(z) = z- \sum_{\nu=k+1}^{\infty} a_{\nu} \lambda_{\nu} z^{\nu}$$
is in class $\;\mathcal{P}_{\tau,\mu}(k,\delta,\xi),\;$
where
\begin{align*}
\xi \leqq  1 - \dfrac { (1 - \gamma)^{2}\;(\mu+1)_{k}} {(1+\delta k ) \, (\tau +1)_{k}}.
\end{align*}
\end{theorem}

\noindent
\textbf{Proof.}
With a view to finding the {\it largest} $\xi$, by supposing that
$\Omega(z) \in  \mathcal{P}_{\tau,\mu}(k,\delta,\xi)$, we have
\begin{align}\label{Inq:1}
\sum_{\nu =k+1}^{\infty} \dfrac{(1+\delta \nu -\delta) \,
\Gamma(\nu+\tau)\Gamma{(\mu+1)}}{( 1 -\xi) \Gamma(\nu+\mu)\Gamma(\tau+1)}
a_{\nu} \lambda_{\nu} \leqq  1 .
\end{align}
Since  $f(z), \psi(z) \in  \mathcal{P}_{\tau,\mu}(k,\delta,\gamma)$,
we know that
\begin{align*}
\sum_{\nu =k+1}^{\infty} \dfrac{(1+\delta \nu - \delta)\,
\Gamma(\nu+\tau)\Gamma{(\mu+1)}}{(1- \gamma)
\Gamma(\nu+\mu)\Gamma(\tau+1)}\; a_{\nu}  \leqq  1
\end{align*}
and
\begin{align*}
\sum_{\nu =k+1}^{\infty} \dfrac{(1+\delta \nu-\delta)\,
\Gamma(\nu+\tau)\Gamma{(\mu+1)}}{(1- \gamma)\Gamma(\nu+\mu)\Gamma(\tau+1)}\;
\lambda_{\nu} \leqq  1.
\end{align*}
Thus, by using the Cauchy-Schwarz inequality, we obtain
\begin{align*}
\sum_{\nu =k+1}^{\infty} \dfrac{(1+\delta \nu - \delta) \,
\Gamma(\nu+\tau)\Gamma{(\mu+1)}}{(1- \gamma)\Gamma(\nu+\mu)\Gamma(\tau+1)}
\;\sqrt{ \lambda_{\nu} a_{\nu} }  \leqq  1,
\end{align*}
which implies that
\begin{align*}
&\dfrac{(1+\delta \nu - \delta) \, \Gamma(\nu+\tau)\Gamma{(\mu+1)}}
{(1- \gamma)\Gamma(\nu+\mu)\Gamma(\tau+1)} \; \sqrt{ \lambda_{\nu} a_{\nu}}\\
&\qquad \quad \leqq  \dfrac{(1+\delta \nu - \delta) \, \Gamma(\nu+\tau)\Gamma{(\mu+1)}}
{( 1 - \xi) \Gamma(\nu+\mu)\Gamma(\tau+1)} a_{\nu} \lambda_{\nu} \leqq 1
\qquad (\nu\geqq k+1).
\end{align*}
that is, that
\begin{align*}
\sqrt{ \lambda_{\nu} a_{\nu} }  \leqq \dfrac{1 -\xi}
{1-\gamma}\qquad (\nu\geqq k+1).
\end{align*}
We note also that
\begin{align*}
\sqrt{\lambda_{\nu} a_{\nu}}
\leqq  \dfrac{(1- \gamma)\Gamma(\nu+\mu)\Gamma(\tau+1)}
{(1+\delta \nu - \delta) \, \Gamma(\nu+\tau)\Gamma(\mu+1)}.
\end{align*}
We now need to show that
\begin{align}
\frac{(1- \gamma)\Gamma(\nu+\mu)\Gamma(\tau+1)}
{(1+\delta \nu - \delta) \, \Gamma(\nu+\tau)\Gamma(\mu+1)}
\leqq \dfrac{1 - \xi}{1 -\gamma}
\end{align}
or, equivalently, that
\begin{align*}
\xi  \leqq 1 -\dfrac { (1 - \gamma)^{2}\;\Gamma(\nu+\mu)\Gamma(\tau+1)}
{(1+\delta \nu - \delta) \, \Gamma(\nu+\tau)\Gamma{(\mu+1)}}.
\end{align*}
Upon letting
\begin{align*}
\Xi(\nu):= 1 -\dfrac{ (1 -\gamma)^{2}\;\Gamma(\nu+\mu)\Gamma(\tau+1)}
{(1+\delta \nu - \delta) \, \Gamma(\nu+\tau)\Gamma{(\mu+1)}},
\end{align*}
we can easily see that the function $\Xi(\nu)$ is non-decreasing in $ \nu$.
We thus obtain
\begin{align*}
\xi \leqq \Xi(k+1) \leqq 1-\dfrac{(1 - \gamma)^{2}\;(\mu+1)_{k}}
{(1+\delta k ) \, (\tau +1)_{k}} \qquad (k\in \mathbb{N}).
\end{align*}
The result is sharp with the extremal function given by
\begin{align}
f(z)= \psi(z)= z-\dfrac{(1 -\gamma)(\mu+1)_{k}}{(1+\delta k ) \,
(\tau +1)_{k}} \; z^{k+1}\qquad (k\in \mathbb{N}).
\end{align}

\section*{Competing interests } The authors declare that they have no competing
interests.

\section*{Author's contributions}
    All the authors jointly worked on deriving the results and approved the final manuscript.

\end{document}